\documentstyle{article}

\title{Computads and slices of operads.}

\author{M.A. Batanin\protect \footnote{The author holds the Scott Russell Johnson Fellowship in
the Centre of Australian Category Theory at Macquarie University}\\ Macquarie University,  North
Ryde, NSW 2109, Australia \\
mbatanin@math.mq.edu.au} 
\date{4 August 2002 }

\newtheorem{theorem}{\bf Theorem}[section]
\newtheorem{defin}{\bf Definition}[section]

\newtheorem{lemma}{\bf Lemma}[section]

\newtheorem{cor}{\bf Corollary}[theorem]

\newcommand{\C}{\mbox{$\cal C$}}
\newcommand{\F}{\mbox{$\cal F$}}
\newcommand{\Proof}{\noindent {\bf Proof. \ }}

\newcommand{\Q}{
{\unitlength=0.25mm
\begin{picture}(500,10)(-10,0)
\put(440,10){\line(0,-1){10}}
\put(440,0){\line(1,0){10}}
\put(450,0){\line(0,1){10}}
\put(450,10){\line(-1,0){10}}
\put(451,11){\line(-1,0){10}}
\put(451,1){\line(0,1){10}}
\put(450.5,0.5){\line(0,1){10}}
\put(450.5,10.5){\line(-1,0){10}}
\end{picture}}}

\newcommand{\Example}{\noindent \makebox[24mm]{{\bf Example
\hspace{-1mm}\addtocounter{example}{1} 
\thesection.\theexample \ }}}

\newtheorem{conj}{\bf Conjecture}[section]

\newcommand{\Remark}{\noindent \makebox[23mm]{{\bf Remark
\hspace{-1mm}\addtocounter{remark}{1} 
\thesection.\theremark \ }}}

\newcommand{\A}{\mbox{\LARGE\it a}}
\newcommand{\B}{
{\unitlength=0.25mm
\begin{picture}(8.5,10)(0,0)

\put(5,6.5){\makebox(0,0){\mbox{\large $b
$}}}
\put(3.2,10.5){\makebox(0,0){\mbox{
$\scriptstyle o
$}}}
\put(4.4,7.5){\makebox(0,0){\mbox{$\scriptstyle\cdot
$}}}
\end{picture}}}

\newcommand{\Ps}{\mbox{$\cal P$}}

\newcommand{\op}{\mbox{\large $\vartheta$}}

\newcommand{\G}{\mbox{$\mathcal G$}}

\newcommand{\Sk}{\mbox{$tr$}}
\newcommand{\U}{\mbox{$\cal W$}}

\newcounter{remark}[section]
\newcounter{example}[section]

\begin{document}

\maketitle

\begin{abstract} For a given $\omega$-operad $A$ on globular sets we 
introduce a sequence of symmetric operads on $Set$ called slices of $A$
and show how the connected limit preserving properties of slices are related to the property of 
the category of $n$-computads of $A$ being a presheaf topos.

\

1991 Math. Subj. Class. 18C20, 18D05      
\end{abstract}

\section{Introduction.}

Computads were invented by Street \cite{StL} as a tool for the presentation of strict
$n$-categories. They attracted a new wave of interest in recent years due to the
development of the theory of weak higher categories. It  also became evident that we
often need some more general types of computads than Street's computads. For example, the
theory of surface diagrams in 3{\small D}-space naturally leads to the use of so called
Gray-computads \cite{MT}.  In our paper \cite{BatP} computads for 
 magma-type globular theories were used. 

In our paper \cite{BatC} we  construct a general theory of computads
for finitary monads on globular sets. An important class of such monads consists of so
called analytic monads \cite{BS} which can be identified with higher operads in $Span$ in
the sense of \cite{BatN}. The examples in the previous paragraph all belong to this class of
monads. 

In \cite{BatC} some properties of computads for analytic monads were established. In
particular, it was claimed that computads form a presheaf topos. This statement in the case
of Street's $2$-computads was proved by Shanuel and then reproved by Carboni and Johnstone
\cite{CJ}. Unfortunately, the proof we gave in \cite{BatC} and \cite{BatP} turned out to be 
incorrect. In \cite{MZ} Makkai and Zawadowski observed that the category of Street's
$3$-computads can not be a presheaf topos.

In this paper we study this question more carefully. We find a sufficient condition when 
computads for a given analytic monad on globular sets do form a presheaf category. The
condition is given in terms of a sequence of symmetric operads  in the category of sets which  we
can construct from the  analytic monad. We call this sequence the sequence of {\it slices of the
operad}. We also show that if the slices are {\it normalised} then the condition is even necessary.

We also give  examples of monads for which this condition is satisfied. A surprising result is
that $n$-computads for weak $n$-categories do form a presheaf category for any $n$. This result
is also true for $3$-computads for Gray-categories.  

It seems to us that the slices of operads are closely related to the coherence problem for weak
$n$-categories and  we suggest a couple of conjectures about it in section \ref{slice}.

\noindent {\bf Acknowledgements.}  I would like to thank Ross Street for
stimulating discussion during my work on this paper. I am also grateful to Michael Makkai and Marek
Zawadowski for informing me about their example, which was a starting point for this work.  
 Finally, I acknowledge the financial support of the Scott Russell Johnson Memorial 
 Foundation and the
Macquarie University Research Commitee.

\section{Computads.}

By an $n$-globular (globular  if $n=\omega$) set we mean a sequence (infinite if $n=\omega$)
of sets 
$$X_{0},X_{1},\ldots,X_{k},\ldots, X_n$$
together with source and target maps
$$s_{r-1},t_{r-1}:X_{r}\longrightarrow X_{r-1}$$
satisfying the equation: 
$$
s_{r-1}\cdot s_{r} = s_{r-1}\cdot t_{r} \ , \
t_{r-1}\cdot s_{r} = t_{r-1}\cdot t_{r}.
$$
The set $X_r$ is called the set of $r$-cells of $X$. Sometimes we will use also notation $(X)_r$
for  this  set. 

Every $(n-1)$-globular set can be considered as an $n$-globular set
with empty set of $n$-cells. 
So we have a chain of inclusion functors
$$ Set = Glob_0\subseteq Glob_1 \subseteq \ldots \subseteq Glob_k \subseteq Glob_{k+1} \ldots \subseteq
Glob$$ and  each of the  inclusion functors 
$$L_k:Glob_k \longrightarrow Glob_n$$ 
has a right adjoint
$$tr_{k}: Glob_{n} \rightarrow Glob_k .$$

Let $A=(A,\mu,\epsilon)$ be a finitary monad on 
 $Glob$. We denote by $A_n$ the $n$-truncation of $A$, i.e. 
the restriction of $A$ to the category $Glob_n$ of $n$-globular sets. 
The category of algebras of
$A_n$ will be denoted by $Alg_n$ and the corresponding forgetful 
functor will be denoted by
$$W_n:Alg_n \longrightarrow Glob_n .$$

We now make  the following inductive definition \cite{BatC}:

\noindent  The category $Comp_0$ of $A_{0}$-computads  is $Glob_0$.  
 The functors $${\cal W}_0 = W_{0}:Alg_0 \rightarrow Comp_0$$ 
 $${\cal F}_{0}= F_{0}:Comp_{0}\rightarrow Alg_{0}$$
 are the forgetful and free $A_{0}$-algebra functors, respectively.

  Let us 
 suppose now that the category $Comp_{n-1}$ of $A_{n-1}$-computads is 
 already defined together with two functors:
 $$ {\cal W}_{n-1}:Alg_{n-1}\rightarrow Comp_{n-1}$$
  $$ {\cal F}_{n-1}:Comp_{n-1}\rightarrow Alg_{n-1}$$
 such that ${\cal F}_{n-1}$ is left adjoint to ${\cal W}_{n-1}$.

  \begin{defin}
  An $A_{n}$-computad $ \cal C$ is a triple $(C,\phi,{ \cal C}')$ 
consisting 
  of an $n$-globular set $C$, an $A_{n-1}$-computad ${ \cal C}'$ and 
an 
  isomorphism 
  $$\phi:W_{n-1}({\cal F}_{n-1}{ \cal C}')\rightarrow {tr}_{n-1}C$$
 in $Glob_{n-1}$.\end{defin}

 Let $G$ be an object of $Alg_{n}$. The counit of the adjunction 
$\F_{n-1}\dashv \U_{n-1}$ gives a morphism
$$r_{n-1}:\F_{n-1}\U_{n-1}\Sk_{n-1}G\rightarrow \Sk_{n-1}G. $$

Define an $n$-globular set $\G$ in the following way. The 
$(n-1)$-skeleton  of $\G$ coincides with 
$W_{n-1}\F_{n-1}\U_{n-1}\Sk_{n-1}G$ and  
$$\G_{n}= \{(\xi,a,\eta) \in \G_{n-1}\times G_{n}\times \G_{n-1} \ | 
\ s_{n-2}\xi = 
s_{n-2}\eta, \ t_{n-2}\xi = t_{n-2}\eta ,\ $$ $$\ \ \ \ \ \ \ \  \ \ 
\ \ \ \ \ \ \ 
s_{n-1}a=r_{n-1}(\xi), \ t_{n-1}a = r_{n-1}(\eta) \}.$$
Define $$s_{n-1}(\xi,a,\eta) = \xi \ , \ t_{n-1}(\xi,a,\eta) = 
\eta . $$
Then put
$$\U_{n}G=(\G,id, \U_{n-1}\Sk_{n-1}G).$$

For an $A_{n}$-computad $\C=(C,\phi,{ 
\C}'), $  define
$$V_{n}(\C) 
= C $$ 
and $V_{0}= id$ for $n=0$. 

 Define a natural transformation
$$\Theta_{n}: V_{n}\U_{n}\rightarrow W_{n},$$
 to be the morphism of $n$-globular sets
which coincides with 
$$W_{n-1}r_{n-1}:W_{n-1}\F_{n-1}\U_{n-1}\Sk_{n-1}G\rightarrow 
W_{n-1}\Sk_{n-1}G$$
up to dimension $n-1$ and has
$$\Theta_{n}(\xi,a,\eta)=a$$
in dimension $n$.

Let us define a new monad $I_A$ on globular sets by means of the following pushout.

\begin{center}{\unitlength=1mm

\begin{picture}(50,30)(5,5)

\put(10,30){\makebox(0,0){\mbox{$L_{n-1}\Sk_{n-1}X$}}}

\put(10,26){\vector(0,-1){12}}

\put(-15,15){\shortstack{\mbox{$ $}}}

\put(10,10){\makebox(0,0){\mbox{$X$}}}

\put(16,10){\vector(1,0){27}}

\put(28,12){\shortstack{\mbox{$ \epsilon_{I_A}$}}}

\put(52,10){\makebox(0,0){\mbox{$I_A X$}}}

\put(52,26){\vector(0,-1){12}}

\put(50,15){\shortstack{\mbox{$ $}}}

\put(52,30){\makebox(0,0){\mbox{$L_{n-1}\Sk_{n-1}A X$}}}

\put(22,30){\vector(1,0){17}}

\put(32,32){\shortstack{\mbox{$ $}}}

\end{picture}}

\end{center}

The algebras of $I_A$ are globular sets together with an $A_{n-1}$-algebra structure on its 
$(n-1)$-truncation. Notice that the categories of $A_n$-computads and $(I_A)_n$-computads are 
canonically isomorphic. Moreover, the functor $V$ together with the $A_{n-1}$-algebra structure 
on $\Sk_{n-1}VC \simeq W_{n-1}({\cal F}_{n-1}{ \cal C}')$ is left adjoint to the forgetful
functor from the category of $I_A$-algebras to $A_n$-computads and $\Theta_n$ is the counit of
this adjunction. So, the functor
${\cal F}_n$ is canonically isomorphic to a composite of $V$ and $\Gamma$ which is  left
adjoint to the
 restriction  functor 
$$l^{\star}: Alg_n \longrightarrow Alg_{I_A}$$
induced by an obvious morphism of monads 
$$l:I_A\rightarrow A_n.$$ 
This left adjoint exists due to the finitary assumption \cite{Kelly}. 

We also can talk about $\omega$-computads.
Recall \cite{BatC} 
that the $n$-truncation of an 
$(n+1)$-computad $(C,\phi,\C)$ is the $n$-computad $\C$. 
\begin{defin} Let $A$ be a finitary monad on $Glob$. An ${\omega}$-computad for $A$
is a sequence
$ \C_n$ of  $n$-computads for $A$ together with a sequence of isomorphisms 
$$c_n: tr_n(\C_{n+1}) \rightarrow \C_n .$$
A morphism of $\omega$-computads is a sequence of morphisms of 
$n$-computads which commutes in the obvious sense with the structure  isomorphisms.  
\end{defin}

We use the techniques of \cite{Kelly} for an explicit  construction of the left adjoint $\Gamma$
into the category of
$A_n$-algebras.

Let $X=M_0$ be an $I_A$-algebra
and let $M_1$ be the following coequalizer in $Glob_n$  

\begin{center}{\unitlength=1mm

\begin{picture}(80,10)(0,15)

\put(69,20){\makebox(0,0){\mbox{$A_nI_A M_{0}$}}}

\put(40,20){\makebox(0,0){\mbox{$ A_n M_{0}$}}}

\put(61,21){\vector(-1,0){15}}

\put(61,19){\vector(-1,0){15}}

\put(54,17){\shortstack{\mbox{$ \eta $}}}

\put(53,22){\shortstack{\mbox{$ Ak $}}}

\put(15,20){\makebox(0,0){\mbox{$  M_{1}$}}}

\put(34,20){\vector(-1,0){15}}

\put(26,21){\shortstack{\mbox{$ \pi_1 $}}}

\end{picture}}

\end{center}

\noindent where $k$ is the $I_A$-algebra structure morphism for $X$ and $\eta$ is the
composite
$\mu\cdot A_n(l)$. Notice, that $k$ is an identity in dimension $n$.

Suppose that the globular set $M_{r},$ together with the morphism 
$$\pi_r: A_n M_{r-1} \rightarrow M_r,$$
are already constructed. Then define $M_{r+1}$ to be 
the following coequalizer.  

\begin{center}{\unitlength=1mm

\begin{picture}(100,25)

\put(99,20){\makebox(0,0){\mbox{$A_n^2 M_{r-1}$}}}

\put(93,16){\vector(-1,-1){9}}

\put(92,10){\shortstack{\mbox{$ \mu_n $}}}

\put(79,4){\makebox(0,0){\mbox{$ A_n M_{r-1} $}}}

\put(76,7){\vector(-1,1){9}}

\put(67,10){\shortstack{\mbox{$\epsilon_n  $}}}

\put(62,20){\makebox(0,0){\mbox{$ A_n^2 M_{r-1}$}}}

\put(91,20){\vector(-1,0){21}}

\put(80,21){\shortstack{\mbox{$ id $}}}

\put(30,20){\makebox(0,0){\mbox{$ A_n M_{r}$}}}

\put(53,20){\vector(-1,0){16}}

\put(40,21){\shortstack{\mbox{$ A_n\pi_r $}}}

\put(3,20){\makebox(0,0){\mbox{$  M_{r+1}$}}}

\put(23,20){\vector(-1,0){15}}

\put(12,21){\shortstack{\mbox{$ \pi_{r+1} $}}}

\end{picture}}

\end{center}

Then we have the following sequence of morphisms
$$M_{0}\stackrel{\epsilon_{n}}{\longrightarrow} A_{n}M_{0}
\stackrel{\pi_{1}}{\longrightarrow} M_{1} 
\stackrel{\epsilon_{n}}{\longrightarrow}
A_{n}M_{1} \stackrel{\pi_{2}}{\longrightarrow} \ldots $$
We denote the colimit of it by $M_{\infty}X$.
According to \cite{Kelly} $M_{\infty}X$ has a natural $A_n$-algebra
 structure
given by   $\pi_{\infty}= \mbox{colim}\hspace{0.5mm} \pi_r$,\ and
this is indeed the  free $A_n$-algebra
generated by $X$.

\section{\label{slice}Suspensions and slices of globular operads.} 

Every strict $\omega$-category has an underlying globular set. 
This functor has a left adjoint
$$D:Glob \longrightarrow \omega\mbox{\it -Cat} .$$
We will also denote by $(D,\mu,\epsilon)$  the monad generated by
this adjunction (notice, that in \cite{BatN} this monad was denoted by $D_s$). In \cite{BatN} a
description of $D$ in terms of plain trees was  presented. 

Recall \cite{StP} that a natural transformation $p:R\rightarrow Q$ between two functors  
is called {\it cartesian} if for every morphism $f:X\rightarrow Y$ the naturality square 

  {\unitlength=1mm

\begin{picture}(60,25)

\put(45,5){\makebox(0,0){\mbox{$Q(X)$}}}

\put(45,16){\vector(0,-1){8}}

\put(47,12){\shortstack{\mbox{$p $}}}

\put(45,20){\makebox(0,0){\mbox{$R(X)$}}}

\put(73,5){\makebox(0,0){\mbox{$Q(Y)$}}}

\put(52,5){\vector(1,0){15}}

\put(56,6){\shortstack{\mbox{$Q(f)$}}}

\put(73,20){\makebox(0,0){\mbox{$R(Y)$}}}

\put(52,20){\vector(1,0){15}}

\put(56,21){\shortstack{\mbox{$ R(f) $}}}

\put(71,16){\vector(0,-1){8}}

\put(73,12){\shortstack{\mbox{$p$}}}

\end{picture}}

\noindent is a pullback. Recall also that an endofunctor $A$ on $Glob$ is called {\it analytic}
 if it is equipped  with
a cartesian natural transformation (augmentation) $p:A\rightarrow D$. Such an endofunctor is 
determined up to isomorphism by  a collection
$$p(1):A(1)\rightarrow D(1),$$
where $1$ is the terminal globular set
and it is connected limits preserving.
A monad on $Glob$ is called analytic if its functor part is analytic and unit and multiplication are
cartesian natural transformation. The category of analytic monads is equivalent to the category of
$\omega$-operads in $Span$. 

The following definition is due to Joyal \cite{J}. 
An endofunctor $\A$ on $Set$ is called {\it analytic} if it can be represented as a `Taylor series' 
$$ \A(X)= \sum_{n\ge 0} A[n]\times_{_{\Sigma_n}} X^n ,$$
where $A[n], n\ge 0$, is a symmetric collection, i.e. a family of sets equipped  with an action of
the symmetric group $\Sigma_n$ on $A[n]$. 
The analytic functors are closed under composition and the monoids in this monoidal categories are
called {\it symmetric operads}. 

Symmetric operads are a special case of algebraic theories in $Set$. Another special case of
algebraic theories called {\it strongly regular} theories was considered by Carboni and Johnstone in
\cite{CJ}. These are theories which can be given by equations without permutations and repetitions
of symbols. For example, the theory of monoids is such a theory, but the theory of commutative
monoids is not. In \cite{CJ} a characterisation of strongly regular theories is established. They
are exactly the theories  given by {\it nonsymmetric operads} in $Set$. The last 
are monoids with respect to composition in the monoidal category of endofunctors of the form
$$\hspace{40mm} \A(X)= \sum_{n\ge 0} A[n]\times X^n , \hspace{40mm} *$$
where $A[n], n\ge 0$ is a nonsymmetric collection, i.e. just a sequence of sets. We will call
the functors of the form $(*)$ {\it strongly analytic}. It was also proved  in \cite{CJ} that
strongly analytic functors preserve connected limits.

\begin{defin} An $n$-globular set $X$ is called $k$-terminal if its 
$k$-th truncation is a terminal $k$-globular set.
An algebra of a monad on $n$-globular sets is called $k$-terminal 
if its underlying globular set is
$k$-terminal.\end{defin}

We denote by $Glob_n^{(k)}$ the category of $k$-terminal $n$-globular sets. 
Clearly, $Glob_n^{(k)}$ is isomorphic to $Glob_{n-k-1}.$ 
For a monad $A$ on $Glob$ we denote by
$Alg_n^{(k)}$ the category of $k$-terminal algebras of $A_n$.  We have a
 restriction of the forgetful  functor $W$
$$W^{(k)}:Alg_n^{(k-1)} \rightarrow Glob_n^{(k-1)},\ k\ge 1
.$$ 
It is not hard to prove that this functor is monadic at least for a 
finitary monad $A$ \cite{W}. 
Hence, we have a monad $S^k A_n $ on $Glob_{n-k}$ such that its category 
of  algebras is equivalent to $Alg_n^{(k-1)}$. We also put $S^0 A = A$.
\begin{defin}\cite{W} $S^k A_n$ is called the $k$-fold suspension of 
$A_n$ \end{defin}
Now if $k=n$ then $S^k A_k $ is a monad on
$Glob_0 = Set$. The proposition 2.1 and the theorem 10.2 from \cite{BatEH} assert that this 
monad is actually  a symmetric operad on $Set$.

\begin{defin}  The symmetric operad $S^k A_k $  will be called the $k$-th slice of
$A$. We will denote this operad 
by $\Ps_k(A)$. 
 \end{defin}

\Example For any operad $A$ its $0$-slice is given by a symmetric operad which
underlying collection consists of a monoid $A[U_0]$ in dimension
$1$ and empty sets in other dimensions. The tree $U_0$ is the only tree of
height $0$.

\

\Example  The first  slice of  the terminal operad $D$ is free monoid operad.
All the higher slices are the free commutative monoid operad.

\

It is proved in \cite[Theorem 10.1]{BatEH} that the first slice of an operad is always  a
free symmetric operad  on some nonsymmetric
operad \cite{BatEH} and, hence, is always a strongly regular theory.

\

\Example For the bicategory operad, its first slice is the nonsymmtric operad freely 
generated by a pointed collection
which has exactly one operation in dimensions $0,1,2$. 
The second slice is the free commutative monoid operad. 

\

\Example  For the Gray operad $G$ \cite{BatN} the first slice is the free monoid operad, 
the second slice is the double-monoid-with-common-unit 
operad i.e. a set with two independent monoid structures 
and  common unit. So $\Ps_2(G)$ is a strongly regular theory.
 The third slice is the free commutative
monoid operad.

\

\Example For a free operad on a globular  collection, the  slices are free
symmetric operads on some nonsymmetric collections and are, therefore, strongly regular
theories. The category of $\omega$-computads for such  operads were used in \cite{BatP}.

\

\Example For the universal contractible $\omega$-operad $K$ from \cite{BatN} the slices are free
symmetric
 operads on  nonsymmetric collections. This can be easily seen from the construction of
$K$ given in \cite{BatN}. Hence, all the slices of $K$ are strongly regular theories.
Recall that the algebras of $K$ are by definition weak $\omega$-categories. 

\

\Example For the universal contractible $n$-operad its slices up to dimension $n-1$  are free
symmetric
 operads on some nonsymmetric collections but its $n$-th slice is the free commutative
monoid operad. The algebras of this operad are weak $n$-categories.

\

In the theory of symmetric operads a very important condition is  freeness of the action of the
symmetric groups. For example,  $E_{\infty}$-operads  are exactly those operads which are
contractible and have free action of the symmetric groups. If the action is not free it usually
means that the corresponding  algebras have some homotopy degeneracy like the vanishing of some
Whithead products or Postnikov invariants.  

From the examples above we see that slices carry with them some information about the homotopy
behaviour of the higher operads. It seems to us that the condition for slices to be regular
theories is the  correct analogue of the  condition of freeness of action of the symmetric
groups. So our conjecture is

\begin{conj} Suppose that  an $n$-operad $A$ is contractible, contains a system of binary
compositions \cite{BatN}, and all its slices up to dimension
$n-1$ are strongly regular theories. Then every weak $n$-category is weakly equivalent to an
$A$-algebra. \end{conj}

At the  time of writing  it is not completely clear what the right notion of 
`semistrict'
$n$-categoryshould be. The desirable properties are: 
\begin{itemize}
\item  every weak $n$-category must be equivalent to a semistrict one;
\item
 the notion of `semistrict' $n$-category is  `minimal'  with the above property.
\end{itemize}

In dimension $2$ this is just the notion of strict $2$-category. In dimension $3$ it is the notion
of Gray-category  \cite{GPS}. Crans has a candidate for dimension $4$ and some ideas about higher
dimensions \cite{Crans}. Here we  risk to suggesting a conjecture.

\begin{conj} There is a unique contractible  $n$-operad $G_n$ with the property that
$\Ps_k(G_n), \ 0\le k \le n-1,$ is the free $k$-fold monoid operad. A semistrict $n$-category is
an algebra for this operad. \end{conj}

\section{Weak limits and coequalisers}

This section has a technical character and contains some elementary facts about weak pullbacks
and coequalisers we will need in next section.

\begin{defin} Let $F:\Lambda\rightarrow C$ be a functor between two categories and let
$W\stackrel{p_{\lambda}}{\rightarrow} F(c_{\lambda})$ be a cone over $F$. It is called a
weak limit of
$F$ if for any other cone $V\stackrel{q_{\lambda}}{\rightarrow}F(c_{\lambda})$ there
exists a morphism
$r:V\rightarrow W$ such that $q_{\lambda}=p_{\lambda}\cdot r$. \end{defin}

\Remark It is obvious that if  limit of a functor $F$ exists then it is a retract of any weak
limit of $F$. Moreover, in order to prove that $W$ is a weak limit it is enough to construct
a section of the canonical morphism from $W$ to the limit of $F$ which makes some obvious
diagrams commutative. We will use  this simple observation extensively.

Following \cite{J} and \cite{W} we call a natural transformation between two functors {\it
weak cartesian} provided every naturality square is a weak pullback. 

\begin{lemma} \label{wcart} Suppose 

{\unitlength=1mm

\begin{picture}(60,12)

\put(28,5){\makebox(0,0){\mbox{$ C$}}}

\put(42,5){\vector(-1,0){10}}

\put(36,6.5){\shortstack{\mbox{\small $p $}}}

\put(45,5){\makebox(0,0){\mbox{$ A$}}}

\put(63,5){\makebox(0,0){\mbox{$B $}}}

\put(59,4.3){\vector(-1,0){10}}

\put(59,5.8){\vector(-1,0){10}}

\put(54,2){\shortstack{\mbox{\small $\chi $}}}

\put(54,7){\shortstack{\mbox{\small $\zeta $}}}

\end{picture}}

\noindent is a coequaliser  of two weakly cartesian transformations between functors
$A,B:\Lambda\rightarrow Set$. Then
$p$ is weakly cartesian. \end{lemma}

\Proof Let $f:X\rightarrow Y$ be a map of sets and let $P$ be the
pullback of $C(f)$ and $p_Y$ i.e.
$$P = \{(c,a) | C(f)(c)= p_{Y}(a) \}. $$
 We have to prove that  there is a section
$s$ of the canonical map
$A(X)\rightarrow P$ which makes the following diagram commutative

{\unitlength=0.9mm

\begin{picture}(60,45)(15,-2)

\put(45,25){\makebox(0,0){\mbox{$ C(X)$}}}
\put(45,20){\vector(0,-1){12}}

\put(87,35){\makebox(0,0){\mbox{$A(X) $}}}
\put(75,25){\makebox(0,0){\mbox{$P $}}}
\put(75,20){\vector(0,-1){12}}
\put(79,34){\vector(-4,-1){25}}
\put(86,32){\vector(-1,-4){5.9}}
\put(77,27){\vector(1,1){5}}

\put(67,25){\vector(-1,0){13}}

\put(35,14){\shortstack{\mbox{\small $C(f) $}}}
\put(83,14){\shortstack{\mbox{\small $A(f) $}}}

\put(45,5){\makebox(0,0){\mbox{$C(Y)$}}}

\put(76,5){\makebox(0,0){\mbox{$A(Y)$}}}

\put(67,5){\vector(-1,0){13}}

\put(63,32){\shortstack{\small \mbox{$p_{\scriptscriptstyle X} $}}}
\put(60,6.5){\shortstack{\small \mbox{$p_{\scriptscriptstyle Y} $}}}

\put(77,29){\shortstack{\mbox{\small $s $}}}

\end{picture}}

Let us take  $(c,a)\in P$ and let $a'\in A(X)$ be such that $p_X(a')= c$.
Put $y=A(f)(a')$.  Then $p_Y(y) = p_Y(a)$. The last equality means
$x$ and $a$ are equivalent with respect to the  equivalence relation generated by
$\chi$ and $\zeta$. Without loss of generality we can assume that there is a
finite sequence
$b_1,\ldots,b_k$ of elements of
$B(Y)$ such that 
$$y=\chi(b_1)\ , \ a = \zeta(b_k), \ \zeta(b_i) = \chi(b_{i+1}).$$  
Since $\chi$ is weakly cartesian we can find a $b'_1\in B(X)$ such that 
$B(f)(b'_1) = b_1$ and $\chi(b'_1)= a'$. Then consider the element
$\zeta(b'_1)$. We have 
$p_X(\zeta(b'_1))= c $ and 
$$ A(f)(\zeta(b'_1)) = \zeta(B(f)(b'_1)) =
\zeta(b_1) = \chi(b_2).$$
Therefore, we can find $b'_2$ such that 
$B(f)(b'_2) = b_2$ and $\chi(b'_2) = \zeta(b'_1).$ Then again 
$p_X(\zeta(b'_2))= c $ and 
$$ A(f)(\zeta(b'_2)) = \zeta(B(f)(b'_2)) =
\zeta(b_2) = \chi(b_3).$$
We can continue this process and finally we get 
$$p_X(\zeta(b'_{k}))= c $$ and 
$$ A(f)(\zeta(b'_{k})) = 
\zeta(b_k) = a.$$
Hence, we can put $s(c,a)= \zeta(b'_k) $. The lemma is therefore proved.

\Q

\begin{lemma}\label{sc} Sequential colimits in $Set$ preserve weak pullbacks.
\end{lemma}

\Proof  It is well known that sequential colimits in $Set$ preserve pullbacks. So it is enough
to prove  that  in a sequential colimit of weak pullbacks we can choose  the
sections of the retractions from pullbacks to weak pullbacks naturally.

 Let us fix
a section $q_i:P_i\rightarrow W_i$ of the canonical 
retraction $W_i\rightarrow P_i$ for every $i\ge 0$. We will construct a
new section
$s_i$ inductively. 

We take $s_0=q_0$. Now
suppose the retraction $s_i$ in the
$i$-th weak pullback

{\unitlength=0.9mm

\begin{picture}(60,45)(15,-2)

\put(50,25){\makebox(0,0){\mbox{$ A_i$}}}
\put(50,22){\vector(0,-1){14}}

\put(83,35){\makebox(0,0){\mbox{$W_i $}}}
\put(72,25){\makebox(0,0){\mbox{$P_i $}}}
\put(72,22){\vector(0,-1){14}}
\put(79,34){\vector(-4,-1){25}}
\put(82,32){\vector(-1,-4){5.9}}
\put(74,27){\vector(1,1){5}}

\put(69,25){\vector(-1,0){15}}

%\put(35,14){\shortstack{\mbox{\small $C(f) $}}}
%\put(83,14){\shortstack{\mbox{\small $A(f) $}}}

\put(50,5){\makebox(0,0){\mbox{$C_i$}}}

\put(73.5,5){\makebox(0,0){\mbox{$B_i$}}}

\put(69,5){\vector(-1,0){15}}

%\put(63,32){\shortstack{\small \mbox{$p_{\scriptscriptstyle X} $}}}
%\put(60,6.5){\shortstack{\small \mbox{$p_{\scriptscriptstyle Y} $}}}

\put(77,28){\shortstack{\mbox{\small $s_i $}}}

\end{picture}}

\noindent is already constructed. Then we can construct $s_{i+1}$ in the
following way. Let $a\in P_{i+1}$ belong to the image of the limiting
map
$\lambda_i:P_i\rightarrow P_{i+1}$ and let us choose a $b\in P_i$ such
that
$\lambda_I (b) = a.$ Then we put $s_{i+1}(a)= w_i(s_i(b)),$ where
$w_i:W_i\rightarrow W_{i+1}$. If $a$ does not belong to $im(\lambda_i)$
then we put $s_{i+1}(a)=q_{i+1}(a).$ 
The sections $s_i$  obviously induce a section
$$\mbox{colim}\hspace{0.3mm} P_i \rightarrow \mbox{colim}\hspace{0.3mm}
W_i$$ of the canonical map $\mbox{colim}\hspace{0.3mm} W_i \rightarrow
\mbox{colim}\hspace{0.3mm} P_i$ which completes the proof. 

\Q

By a similar diagram-chase method one can easily prove the following 
lemma.

\begin{lemma}\label{mono} Suppose  that in a commutative diagram of
coequalisers

{\unitlength=1mm

\begin{picture}(60,29)

\put(28,5){\makebox(0,0){\mbox{$ C_1$}}}

\put(42,5){\vector(-1,0){10}}

\put(36,6.5){\shortstack{\mbox{\small $p_1 $}}}

\put(45,5){\makebox(0,0){\mbox{$ A_1$}}}

\put(63,5){\makebox(0,0){\mbox{$B_1 $}}}

\put(59,4.3){\vector(-1,0){10}}

\put(59,5.8){\vector(-1,0){10}}

%\put(54,2){\shortstack{\mbox{\small $\chi $}}}

\put(25,11.5){\shortstack{\mbox{\small $\phi $}}}

\put(28,20){\makebox(0,0){\mbox{$ C_2$}}}

\put(42,20){\vector(-1,0){10}}

\put(36,21.5){\shortstack{\mbox{\small $p_2 $}}}

\put(60,11.5){\shortstack{\mbox{\small $\psi $}}}
\put(42.5,11.5){\shortstack{\mbox{\small $\zeta $}}}

\put(45,20){\makebox(0,0){\mbox{$ A_2$}}}

\put(63,20){\makebox(0,0){\mbox{$B_2 $}}}

\put(59,19.3){\vector(-1,0){10}}

\put(59,20.8){\vector(-1,0){10}}

%\put(54,2){\shortstack{\mbox{\small $\chi $}}}

%\put(54,7){\shortstack{\mbox{\small $\zeta $}}}

\put(28,16){\vector(0,-1){8}}
\put(45,16){\vector(0,-1){8}}
\put(63,16){\vector(0,-1){8}}

\end{picture}}

\noindent both right commutative squares are weak pullbacks and $\psi$ and
$\zeta$ are monomorphisms, then colimiting map $\phi$ is a monomorphism.
\end{lemma} 

\Q

The following lemma is obvious.
\begin{lemma} If a commutative square  is weakly cartesian and
one of the limiting maps is a monomorphism then the square is cartesian.  \end{lemma}

\Q

\begin{lemma}  Analytic functors on $Glob$ preserve connected weak
limits.\end{lemma}

\Proof  Let $A$ be an analytic functor on $Glob$ and let $C$ be a weak connected  limit of a
diagram of globular sets
$F:\Lambda \rightarrow Glob.$  Then there is a retraction 
$$r:C\rightarrow \lim_{\Lambda}F.$$ 
Hence, we have a retraction 
$$A(r):A(C)\rightarrow A(\lim_{\Lambda}F) \simeq \lim_{\Lambda}A(F) $$
which proves the lemma. 

\Q

\begin{lemma}\label{mixture} Let $\phi: A \rightarrow B$ be a natural transformation in
$Glob_n$ such that $tr_{n-1}\phi$ is cartesian and $(\phi)_n:(A)_n \rightarrow
(B)_n$ is weakly cartesian in $Set$. Then $\phi$ is weakly cartesian.
\end{lemma}

\Proof Let $f:X\rightarrow Y$ be a morphism of globular sets and let $P$ be a 
pullback of $\phi$ and $B(f)$. Then we can assume that $tr_{n-1}P = tr_{n-1}A$.
 Let $\psi: (P)_n \rightarrow (A(X))_n$ be a section of the canonical retraction
$(r)_n:(A(X))_n\rightarrow (P)_n$ which exists due to the weak cartesianness of $(\phi)_n$.
We have to prove that $\psi$ respects source and target operators. 

Indeed, consider a map $\alpha= s_{n-1}(\psi):(P)_n\rightarrow (A(X))_{n-1}$. Then we 
have
$$ s_{n-1}(p_A)_n = s_{n-1}(p_A)_n((r)_n(\psi)) =
s_{n-1}((Af)_n(\psi))= (Af)_{n-1}(\alpha).$$
Analogously 
$$s_{n-1}(p_B)_n = (\phi)_{n-1}(\alpha) ,$$
where $p_A, p_B$ are canonical projections from the pullback $P$. Since 
$(Af)_{n-1}$ and $(\phi)_{n-1}$ are also projections of a pullback we have by its
universal property that $\alpha$ must coincide with $s_{n-1}:(P)_n\rightarrow
(A(X))_{n-1}$. So $\psi$ commutes with the source operator. Analogously it commutes
with target operator. 

\Q

\begin{lemma}\label{map1} Suppose $\phi: \A\rightarrow \B$ is a weakly
cartesian  transformation between two strongly analytic functors  in
$Set$, then it is cartesian. \end{lemma}
\Proof By \cite{CJ} and a theorem of Joyal \cite{J,W} we can assume that $\A$ and $\B$ both are
given by free symmetric collections  $A[n] = \alpha[n]\times \Sigma_n$, $B[n]=
\beta[n]\times \Sigma_n$ and  $\phi$ is given by  equivariant maps of  symmetric collections
$$\phi[n]:\alpha[n]\times \Sigma_n \rightarrow \beta[n]\times \Sigma_n \ , n \ge 0 \ .$$
The map $\phi[n]$ is determined in its turn by a map of nonsymmetric collections
$$\psi[n]:\alpha[n] \rightarrow \beta[n]\times \Sigma_n.$$
Then the natural transformation $\phi$  is the coproduct over $n$ of the composites
$$(\alpha[n]\times \Sigma_n)\times_{_{\Sigma_n}} X^n
\simeq \alpha[n]\times X^n
\stackrel{\psi[n]\times 1}{-\!\!\!\longrightarrow} 
\beta[n]\times \Sigma_n \times X^n \stackrel{1\times k}{\rightarrow}$$
 $$\stackrel{1\times k}{\rightarrow} \beta[n]\times X^n \simeq
(\beta[n]\times
\Sigma_n)\times_{_{\Sigma_n}} X^n , \hspace{30mm} $$
where $k$ is the action of $\Sigma_n$ on $X^n$. Then for the unique map
$X\rightarrow 1$ we have the following commutative naturality diagram

{\unitlength=1mm

\begin{picture}(60,29)(-10,0)

\put(7,5){\makebox(0,0){\mbox{$ A[n]$}}}

\put(15,5){\vector(1,0){15}}

\put(20,6.5){\shortstack{\mbox{\small $\psi[n]$}}}

\put(45,5){\makebox(0,0){\mbox{$B[n]\times \Sigma_n$}}}

\put(83,5){\makebox(0,0){\mbox{$B[n] $}}}

\put(59,5){\vector(1,0){13}}

%\put(54,2){\shortstack{\mbox{\small $\chi $}}}

\put(25,11.5){\shortstack{\mbox{\small $ $}}}

\put(7,20){\makebox(0,0){\mbox{$A[n]\times X^n$}}}

\put(15,20){\vector(1,0){15}}

\put(17,21.5){\shortstack{\mbox{\small $\psi[n]\times 1 $}}}

\put(45,20){\makebox(0,0){\mbox{$ B[n]\times \Sigma_n\times X^n$}}}

\put(83,20){\makebox(0,0){\mbox{$B[n]\times X^n $}}}

\put(59,20){\vector(1,0){13}}
\put(62,21.5){\shortstack{\mbox{\small $1\times k $}}}

%\put(54,2){\shortstack{\mbox{\small $\chi $}}}

%\put(54,7){\shortstack{\mbox{\small $\zeta $}}}

\put(7,16){\vector(0,-1){8}}
\put(45,16){\vector(0,-1){8}}
\put(83,16){\vector(0,-1){8}}

\end{picture}}

\noindent In this diagram both left and right squares are obviously
pullbacks; hence, so is the big square. This is enough to imply the transformation
$\phi$ is cartesian. 

\Q

\begin{lemma} \label{mon2} Let $A$ be an analytic functor on $Glob$ and
let
$f:X\rightarrow Y$ be a map of globular sets such that for a fixed $n\ge 0$ the map $(f)_n$ is a
monomorphism. Then $(A(f))_n$ is a monomorphism. \end{lemma}
\Proof Since $A$ is strongly analytic it is sufficient to prove the lemma for
the case $A=D$. Then it is obvious from the construction of $D$ given in
\cite{BatN}.

\Q

  \section{Computads and slices of operads.}

Let $$k_n: A_n(\F_{n}) \rightarrow \F_{n}$$
be a natural transformation which is given on a computad $\C$ by the
 structure map of 
the algebra $\F_n(\C)$. 

\begin{theorem} Suppose for an $n$-operad $A$ all 
the slices $\Ps_k(A), \ 0\le k \le n, $ are 
strongly regular theories, then 
\begin{itemize}
\item $k_n$ is a cartesian natural transformation;
\item $\F_n$ preserves  connected limits. \end{itemize}
 \end{theorem}

\Proof We will prove the theorem by induction. If $n=0$ the proposition
 is obvious because the
 $0$-operads are just monoids and $0$-computads are sets.  

We assume, therefore, that the natural transformation 
$$tr_{n-1}k_n = k_{n-1}: A_{n-1}\F_{n-1} \rightarrow \F_{n-1}$$ 
is cartesian and $\F_{n-1}$ preserves connected limits.

Now we can use Kelly's method to construct the left adjoint  to the
 restriction
functor 
$$Alg_n^{(n-1)} \rightarrow Alg^{(n-1)}_{I_A} $$

First of all  observe that for the operad $I_A$ the 
natural transformation   
   $$\kappa:I_A (-) \rightarrow (-)$$
is cartesian on the category of $(n-1)$-terminal $I_A$-algebras because
$tr_{n-1}\kappa $ is the constant map
$$1:A_{n-1}(1)\rightarrow 1$$
 and 
$\kappa$ is an identity in dimension $n$.
Hence, $A(\kappa)$ is cartesian. 

By lemma \ref{wcart} the resulting colimit map 
$$\pi_1:A_n M_0 \rightarrow M_1$$
is weakly cartesian in dimension n and,
therefore, by lemma \ref{mixture} $\pi_1$ is weakly cartesian 
because $tr_{n-1}\pi_1 = 1 .$ 

Analogously we have that in  Kelly's construction all 
$\pi_r$ are weakly cartesian transformations
and all $M_r$ preserve connected limits. 

 The last sequential colimit of Kelly's construction 
$$\mbox{colim}\hspace{0.5mm} \pi_r:A_n M_{\infty}\rightarrow M_{\infty}$$
is weakly cartesian by lemmas \ref{mixture} and \ref{sc}. This map is
obviously the map 
$$j : A_n(\Ps_{n}) \longrightarrow \Ps_{n}$$
given on $X$ by the structure morphism of the algebra $\Ps_{n}(X)$.
 Since $\Ps_{n}$ is strongly reguilar theory the
functor
$(A_n(P_{n}))_n$  is strongly
analytic as well. Hence, by lemma \ref{map1} $j$ is even cartesian.

Now let
$$\widehat{(-)}: Alg_{I_A}\rightarrow Alg^{(n-1)}_{I_A}$$
be a functor which assigns to an $I_A$-algebra $X$ the $(n-1)$-terminal
$I_A$-algebra
$\widehat{X}$ with $(\widehat{X})_n = (X)_n$. We obviously have a natural
morphism of $I_A$-algebras $X\rightarrow \widehat{X}$. 

For a computad $\C$ we therefore have a coequalisers diagram

{\unitlength=1mm

\begin{picture}(60,29)

\put(25,5){\makebox(0,0){\mbox{$ N_1$}}}

\put(38,4.8){\vector(-1,0){10}}

%\put(36,6.5){\shortstack{\mbox{\small $p_1 $}}}

\put(45,5){\makebox(0,0){\mbox{$ A_n(\widehat{V\C})$}}}

\put(70,5){\makebox(0,0){\mbox{$A_n I_A(\widehat{V\C}) $}}}

\put(61.5,3.9){\vector(-1,0){10}}

\put(61.5,5.4){\vector(-1,0){10}}

%\put(54,2){\shortstack{\mbox{\small $\chi $}}}

%\put(25,11.5){\shortstack{\mbox{\small $\phi $}}}

\put(25,20){\makebox(0,0){\mbox{$ M_1$}}}

\put(38,20){\vector(-1,0){10}}

%\put(36,21.5){\shortstack{\mbox{\small $p_2 $}}}

%\put(60,11.5){\shortstack{\mbox{\small $\psi $}}}
%\put(42.5,11.5){\shortstack{\mbox{\small $\zeta $}}}

\put(45,20){\makebox(0,0){\mbox{$  A_n({V\C})$}}}

\put(70,20){\makebox(0,0){\mbox{$A_n I_A({V\C}) $}}}

\put(61.5,19.3){\vector(-1,0){10}}

\put(61.5,20.8){\vector(-1,0){10}}

%\put(54,2){\shortstack{\mbox{\small $\chi $}}}

%\put(54,7){\shortstack{\mbox{\small $\zeta $}}}

\put(25,16){\vector(0,-1){8}}
\put(45,16){\vector(0,-1){8}}
\put(70,16){\vector(0,-1){8}}

\end{picture}}
  
\noindent In this diagram the two right vertical morphisms are monomorphisms in
dimension $n$ by lemma \ref{mon2} and, therefore, the colimiting map
is a monomorphism in dimension $n$ by lemma \ref{mono}. In addition, the left
square is a weak pullback by lemma \ref{wcart}. 

What we have here is a map of the first stages of the Kelly machine for 
$V\C$ and $\widehat{V\C}$ respectively. Continuing this process we have as an
output of the Kelly machine in  dimension $n$, a weak pullback 

{\unitlength=1mm

\begin{picture}(60,29)(-13,0)

\put(25,5){\makebox(0,0){\mbox{$(\Ps_{n}\widehat{V\C})_n$}}}

\put(44,4.8){\vector(-1,0){11}}

\put(55,5){\makebox(0,0){\mbox{$( A_n(\Ps_{n}\widehat{V\C}))_n$}}}

\put(25,20){\makebox(0,0){\mbox{$ (\F_n\C)_n$}}}

\put(45,20){\vector(-1,0){14}}

\put(55,20){\makebox(0,0){\mbox{$ ( A_n(\F_n\C))_n$}}}

\put(25,16){\vector(0,-1){8}}
\put(55,16){\vector(0,-1){8}}

\end{picture}}
  
\noindent with vertical morphisms being monomorphisms. So it is a pullback.
By a similar argument, the natural transformation
$$(\F_n\C)_n \longrightarrow (A_n(\Ps_{n}\widehat{V\C}))_n$$
is cartesian.

For a computad morphism $f:\C\rightarrow \C'$ we have now the  following
commutative cube.

{\unitlength=1mm

\begin{picture}(60,50)(-10,-4)

\put(10,25){\makebox(0,0){\mbox{$( \F_n(\C))_n$}}}
\put(10,22){\vector(0,-1){14}}
\put(12,15){\shortstack{\mbox{$ $}}}

\put(35,25){\vector(-1,0){17}}

\put(23,26){\shortstack{\mbox{$ $}}}

\put(45,25){\makebox(0,0){\mbox{$ (A_n(\F_n\C))_n$}}}
\put(45,22){\vector(0,-1){14}}

\put(57,21){\shortstack{\mbox{$ $}}}

\put(57,28){\shortstack{\mbox{\small $ $}}}

\put(10,5){\makebox(0,0){\mbox{$ (\Ps_{n}(\widehat{V\C}))_n$}}}

\put(33,5){\vector(-1,0){14}}

\put(23,6){\shortstack{\mbox{$ $}}}

\put(45,5){\makebox(0,0){\mbox{$ (A_n(\Ps_{n}\widehat{V\C}))_n$}}}

\put(57,1){\shortstack{\mbox{$ $}}}

\put(57,8){\shortstack{\mbox{\small $ $}}}

\put(-3,3){\begin{picture}(50,30)

%back square------------------------------------------

\put(30,35){\makebox(0,0){\mbox{$ (\F_n(\C'))_n$}}}

\put(30,32){\line(0,-1){9.5}}
\put(30,21.2){\vector(0,-1){3}}

\put(32,25){\shortstack{\mbox{$ $}}}

\put(54,35){\vector(-1,0){16}}

\put(43,36){\shortstack{\mbox{$ $}}}

\put(65,35){\makebox(0,0){\mbox{$(A_n(\F_n\C'))_n$}}}
\put(65,32){\vector(0,-1){14}}

\put(77,31){\shortstack{\mbox{$ $}}}

\put(77,38){\shortstack{\mbox{\small $ $}}}

\put(30,15){\makebox(0,0){\mbox{$(\Ps_{n}(\widehat{V\C'}))_n$}}}

\put(53,15){\line(-1,0){4.4}}
\put(47.3,15){\vector(-1,0){6.4}}

\put(43,16){\shortstack{\mbox{$ $}}}

\put(65,15){\makebox(0,0){\mbox{$  (A_n(\Ps_{n}\widehat{V\C'}))_n$}}}

\put(77,11){\shortstack{\mbox{$ $}}}

\put(77,18){\shortstack{\mbox{\small $ $}}}

\end{picture}}

\put(13,28){\vector(1,1){7}}

\put(50,28){\vector(1,1){7}}

\put(13,8){\vector(1,1){7}}

\put(50,8){\vector(1,1){7}}

\end{picture}}

In this diagram the front and rear vertical squares are  pullbacks.
 The bottom horizontal square is a
pullback because $j_n$ is cartesian.
Hence, we have that the top horizontal square is a
pullback in dimension n.  It is also a pullback after truncation by our
inductive assumption. So we have proved that $k_n$ is cartesian.  

Finally, we have to prove that $\F_n$ preserves connected limits. To do this it
is sufficient to prove this result in  dimension $n$.

Let $\C$ be a connected limit of computads $\C_{\lambda}$ and let
$c_{\lambda}:\C\rightarrow \C_{\lambda}$  be the canonical projection.
So we have a  cartesian square

{\unitlength=1mm

\begin{picture}(60,29)(-13,0)

\put(25,5){\makebox(0,0){\mbox{$(\Ps_{n}\widehat{V\C_{\lambda}})_n$}}}

\put(47,4.8){\vector(-1,0){13}}

\put(55,5){\makebox(0,0){\mbox{$( \Ps_{n}\widehat{V\C})_n$}}}

\put(25,20){\makebox(0,0){\mbox{$ (\F_n\C_{\lambda})_n$}}}

\put(48,20){\vector(-1,0){15}}

\put(55,20){\makebox(0,0){\mbox{$ ( \F_n\C)_n$}}}

\put(25,16){\vector(0,-1){8}}
\put(55,16){\vector(0,-1){8}}

\end{picture}}

\noindent But $( \Ps_{n}\widehat{V\C})_n$ is naturally isomorphic to 
$( \Ps_{n}\lim(\widehat{V\C_{\lambda}}))_n $ because  $V$ obviously
preserves limits in dimension $n$. So, after the limit we have a pullback

{\unitlength=1mm

\begin{picture}(60,29)(-15,0)

\put(22,5){\makebox(0,0){\mbox{$\lim((\Ps_{n}\widehat{V\C_{\lambda}}))_n$}}}

\put(45,4.8){\vector(-1,0){9}}

\put(55,5){\makebox(0,0){\mbox{$( \Ps_{n}\widehat{V\C})_n$}}}

\put(22,20){\makebox(0,0){\mbox{$\lim ((\F_n\C_{\lambda}))_n$}}}

\put(48,20){\vector(-1,0){15}}

\put(55,20){\makebox(0,0){\mbox{$ ( \F_n\C)_n$}}}

\put(22,16){\vector(0,-1){8}}
\put(55,16){\vector(0,-1){8}}

\end{picture}}

\noindent where the bottom arrow is an isomorphism because $\Ps_{n}$ preserves
connected limits. So  the top arrow is, and we completed the proof of the
theorem.

\Q

\begin{theorem} Suppose that for an operad $A$  the slices $\Ps_{k}(A), 0\le
k\le n-1$ are strongly regular theories. Then the category of $n$-computads of
$A$ is a presheaf topos. \end{theorem}

\Proof The proof generalizes  example 3.6 from \cite{CJ}. 
We use  induction on $n$. If $n=0$ the statement is true by 
definition. 

Suppose we know that the category $Comp_{n-1}$ is  a 
presheaf topos. Consider the 
 functor 
$$T_{n-1}: Comp_{n-1}\longrightarrow Set \ ,$$
which assigns to a computad $\C$ the set of  parallel pairs of 
$(n-1)$-cells from $W_{n-1}\F_{n-1}\C$. Then we have the equivalence 
of  categories
$$Comp_{n} \sim Set\downarrow T_{n-1}.$$

 Now we prove 
that 
$T_{n-1}$ preserves connected limits. 
Notice 
that $T_{n-1}$ is isomorphic to the following composite
$$ Comp_{n-1}\stackrel{\scriptscriptstyle
{\cal{F}}_{n-1}}{-\!\!\!\longrightarrow} Alg_{n-1} 
\stackrel{Alg_{n-1}(A_{n-1}S^{n-1},-)}{-\!\!\!-\!\!\!-\!\!\!-\!\!\!-\!\!\!-\!\!\!-
\!\!\!-\!\!\!-\!\!\!-\!\!\!-\!\!\!\longrightarrow}Set$$
where $S^{n-1}$ is the  $n-1$-globular set which has two elements $-$ 
and $+$ in 
every dimension and $$s(-) = s(+) = - \ , \ t(-) = t(+) = + . $$
 By our  assumption, $\F_{n-1}$ preserves pullbacks (wide 
pullbacks), so $T_{n-1}$ does. 
 According to the results of \cite{CJ} this is sufficient for
  $Set\downarrow T_{n-1}$ to be a preasheaf topos.

\Q

\begin{cor} The following categories of computads are presheaf toposes:
\begin{itemize}
\item the category of Street $2$-computads (Shanuel, Carboni-Johnstone
\cite{CJ});
\item the category of Gray-computads \cite{MT} and the category of $3$-computads
for
 Gray-categories;
\item the category of $k$-computads for weak $n$-categories for all $0\le k\le
n$;
\item the category of $k$-computads for $P$-magmas \cite{BatP}.
\end{itemize}
\end{cor}

\Proof See examples in section \ref{slice}. 

\Q

The following theorem extends the example of Makkai-Zawadowski.

\begin{theorem} Let $A$ be an operad  such that its slices $\Ps_k(A),\
1\le k\le n-2,$ are normalised in the sense that $\Ps_k(A)[0]=1$. Then
the category of
$n$-computads is a presheaf topos if and only if all the slices
$\Ps_k(A),\ 0\le k \le n-1,$ are strongly regular theories. \end{theorem}

\Proof We only need to prove the only if part of the theorem. For this we will
show that if there exists $\Ps_k(A)$ which is not strongly regular then the
category of $k$-computads is not a presheaf topos; this implies that 
the category of $n$-computads is not a presheaf topos either. So without loss
of generality we can assume that  the category of $(n-1)$-computads is
a presheaf topos but $\Ps_{n-1}(A)$ is not strongly regular, in particular, it is
not connected limits preserving.

Let $\op_{k}, \ 0\le k \le n-2 $ be a  $k$-computad defined by induction
$$\op_{k}= (O_{k-1},id,\op_{k-1})$$
and $\op_0 = 1$, where $O_{k-1}$ is a $(k-1)$-terminal $k$-globular set with 
empty set of cells of  dimension $k$. For this definition
to be valuable, we have to prove that
$$\F_{k}\op_{k} = 1.$$
If $k=0$ it follows from $A_0(1)= 1$. 
Suppose we have proved it up to dimension
$k-1$. Then by applying the Kelly machine we see that the calculation of
$\F_{k}\op_{k}$ amounts to the calculation of a free $A_k$-algebra on the
$I_{A_k}$-algebra $O_{k-1}$. So  the algebra
$\F_{k}\op_{k}$ is isomorphic to $\Ps_k(A)(\emptyset)=1$.

Let us consider the full subcategory of $Comp_{n-1}$ consisting of computads $\C$
with $tr_{n-2}\C= \op_{n-2}.$ Obviously, this subcategory is isomorphic to the
category of sets. By the above argument, the restriction of $\F_{n-1}$ to this 
subcategory is isomorphic to $\Ps_{n-1}(A)$ and, hence, is not connected limit
preserving. So the functor $T_{n-1}$ is not connected limit preserving either
and hence $Comp_{n} \sim Set\downarrow T_{n-1}$ can not be a presheaf topos
 by a theorem from \cite{CJ} again. 

\Q

\begin{cor}[Makkai-Zawadowski \cite{MZ}] 
The category of Street \\ $n$-computads for $n\ge 3$ is not a presheaf topos.
\end{cor}

\end{document}